\documentclass[11pt,leqno]{amsart}

\usepackage{fullpage}
\usepackage{amssymb}
\usepackage{amsmath}
\usepackage{amsxtra}
\usepackage{amscd}
\usepackage{graphicx}
\usepackage{amsfonts}
\usepackage{pb-diagram}

\usepackage{float}
\usepackage{enumerate}
\usepackage{dsfont}
\usepackage{multirow}
\usepackage{color}
\usepackage{rotating}
\usepackage{url}
\usepackage{hyperref}
\usepackage[utf8]{inputenc}
\usepackage[all]{xy}
\usepackage{tikz}
\tikzset{>=latex}
\usetikzlibrary{arrows,decorations.markings,positioning}
\usepackage{tikzinclude}
\usepackage{braids}

\usepackage{wrapfig}

 \usepackage[british,UKenglish,USenglish,american]{babel}
\usepackage{bicaption}

\makeatletter
\let\@wraptoccontribs\wraptoccontribs
\makeatother
\definecolor{darkgreen}{rgb}{0,0.4,0.1}
\hypersetup{
   colorlinks=true,       	% false: boxed links; true: colored links
   linkcolor=red,          % color of internal links (change box color with linkbordercolor)
   citecolor=darkgreen,    % color of links to bibliography
   filecolor=magenta,      % color of file links
   urlcolor=blue			% color of external links
}

\numberwithin{equation}{section}
\newtheorem{thm}{Theorem}

\begin{document}

\title[Braid representations of knots and links]{Diagrammatic representations of knots and links as closed braids}

\author{Sofia Lambropoulou}
\address{ Department of Mathematics,
National Technical University of Athens,
Zografou Campus, GR-157 80 Athens, Greece.}
\email{sofia@math.ntua.gr}
\urladdr{http://www.math.ntua.gr/$\tilde{~}$sofia}

\date{}
\maketitle

\section{Introduction}

% A knot is a homeomorphic image of the circle in 3-space and a link on $k$ components is a homeomorphic image of $k$ circles in 3-space. By `homeomorphic' is meant `bijective and continuous, with its inverse also continuous'. As a knot is just a link on one component, we can say `links' instead of `knots and links'.   The complete classification of knots and links is one of the big open problems in Topology, involving the full complexity of  3-dimensional space. `Classification' means a complete list of non-equivalent such objects. We recall that two knots or links are `equivalent' or `ambient isotopic' or just `isotopic' if one can be transformed into the other via a continuous `elastic'  deformation of the ambient 3-space. So, by `link' we mean the whole of its equivalence class.

\begin{wrapfigure}{R}{0.30\textwidth} 
\centering 
\includegraphics[width=0.25\textwidth]{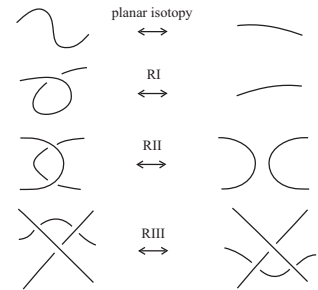} 
\caption{The  Reidemeister moves} 
\label{reidem} 
\end{wrapfigure}

 A way to study knots and links is by studying their regular projections on a plane, called `diagrams'. One big advancement in Low-dimensional Topology was the `discretization' of link isotopy  through the  well-known moves on knot and link diagrams discovered by K. Reidemeister in 1927 \cite{Rd1}, see Fig.~\ref{reidem}.

 The crucial one is the Reidemeister III move, which in terms of dual planar graphs corresponds to the so-called `star-triangle relation'. A second advancement was the `algebraization' of link isotopy, by representing knots and links via braids. Braids are both geometric-topological and algebraic objects; geometric as sets of interwinding, non-intersecting descending strands, algebraic as elements of the so-called Artin braid groups \cite{Ar1,Ar2}. In Fig.~\ref{closures} middle, we can see an example of a braid on 3 strands. The aim of this article is to detail on the connection between links and braids, which is marked by two fundamental results: the {\it Alexander theorem}, whereby links are represented via braids, and the {\it Markov theorem}, which provides the braid equivalence that reflects link isotopy. See Theorems~\ref{alex} and~\ref{markov}.  In 1984, these theorems played a key role in the discovery of a novel link invariant, the Jones polynomial \cite{Jo1,Jo2}.

\section{From braids  to links}

\begin{wrapfigure}{L}{0.6\textwidth} 
\centering 
\includegraphics[width=0.58\textwidth]{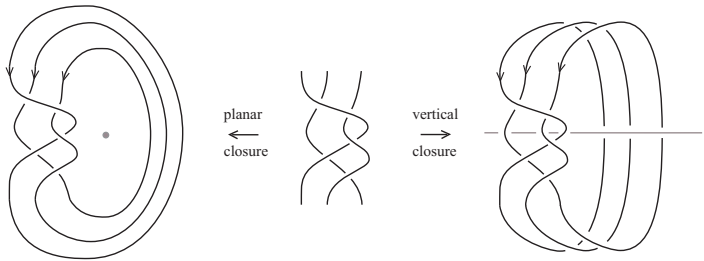} 
\caption{A braid on 3 strands and two closures}
\label{closures} 
\end{wrapfigure} 

More rigorously, a classical  braid on $n$ strands is, geometrically, a homeomorphic image of $n$ arcs in the interior of a thickened rectangle starting from $n$  top points and running monotonically down to $n$  corresponding points, that is, with no local maxima or minima. Each one of the two sets of  endpoints can be assumed colinear  and both of them  coplanar. For an example see middle illustration of  Fig.~\ref{closures}. So a braid can be identified with the braid diagram on the projection plane of its endpoints. Moreover, two braids are {\it isotopic} through any isotopy move of its arcs that preserves the braid structure. These isotopies comprise the Reidemeister moves II and III for braids and planar isotopies, which include the change of relative heights for two non-adjacent crossings, as well as small shifts of endpoints that preserve their order. The  Reidemeister I move cannot apply as the kink introduces local maxima and minima. In the set of braids isotopic elements are considered equal.                                                                            

 The interaction between braids and links takes place through their diagrams. An operation that connects braids and links is the `closure' operation. It is realized by joining with simple non-interwinding arcs the corresponding pairs of endpoints of the braid.  The result is an oriented link diagram that winds around a specified axis, the {\it braid axis}, in the same sense, say the counterclockwise. The link orientation is induced by the top-to-bottom direction of the braid.  There are many ways for specifying the braid axis.  Fig.~\ref{closures} illustrates two of them: in the left-hand illustration the braid axis can be considered to be perpendicular to the braid projection plane, piercing it at one point, its `trace'; in the right-hand illustration the braid axis is a horizontal line parallel to and behind the projection plane of the braid. The two closures are the {\it planar closure} and the {\it vertical closure} respectively.   Clearly, any two  closures of the same braid are isotopic.

\section{From links to braids}

 Now  the following question arises naturally: can one always do the converse? That is, given an oriented  link diagram can one turn it into an isotopic closed braid?  

\begin{wrapfigure}{L}{0.45\textwidth} 
\centering 
\includegraphics[width=0.4\textwidth]{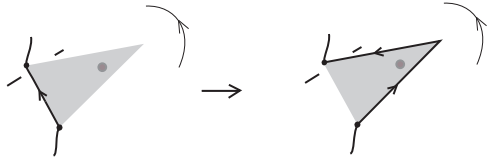} 
\caption{Alexander's braiding of an opposite arc}
\label{tooth} 
\end{wrapfigure}

Note that our diagram may already be in closed braid form for some choice of braid axis. Yet, with a different specified axis it may not be braided any more. Imagine in Fig.~\ref{closures} the axis trace to be placed in some other region of the plane. The answer to the above question is `yes' and the idea is quite simple. Indeed, we first specify a point on the  plane of the diagram, the trace of the braid axis, and we define a `good' direction around this point.  
We then subdivide the arcs of the link diagram into smaller arcs, by marking the transition from `good' arcs to `opposite' arcs, according to whether they agree or not with the  good direction. The  good arcs are left alone. 

\begin{wrapfigure}{R}{0.23\textwidth} 
\centering 
\includegraphics[width=0.18\textwidth]{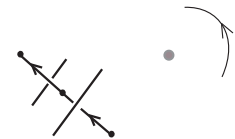} 
\caption{Subdividing an opposite arc}
\label{subdivision} 
\end{wrapfigure}

\noindent An opposite arc is to be subdivided further if needed (Fig.~\ref{subdivision}), so that each subarc can be part of the boundary of a  `sliding triangle'  across the axis trace. See Fig.~\ref{tooth}. The sliding triangle of a subarc does not intersect any other arcs of the link diagram. If the subarc lies over other arcs of the diagram, so does its triangle. Similarly, if it lies under other arcs then its triangle also lies under these arcs. (The notion of the sliding triangle was first introduced by Reidemeister \cite{Rd1}  in the general context of his study of isotopy, not of braiding, and he called these isotopy generating moves {\it $\Delta$-moves}.)
 So, by an isotopy move we can replace the opposite subarc by two good arcs, the other two sides of the sliding triangle, and the opposite subarc is eliminated. The opposite subarcs are finitely many, so the process terminates with a closed braid. The  algorithm outlined above is J.W. Alexander's proof of his homonymous theorem \cite{Al}:

\begin{thm}[J.W. Alexander, 1923] \label{alex}  
Every oriented link diagram may be isotoped to a closed braid.
\end{thm}

Yet, the idea of turning a knot into a braid goes back to Brunn in 1897 \cite{Br}, who observed that any knot has a projection with a single multiple point; from this follows immediately (by appropriate small perturbations) that the diagram can be brought to closed braided form.  In 1974 Joan Birman made  Alexander's braiding algorithm more technical \cite{Bi1} with the purpose of providing a  rigorous proof of the Markov theorem (see Section~\ref{braidequiv}). 

 Another  proof of the Alexander theorem, very appealing to the imagination, is the one by Hugh Morton (\cite{Mo}, 1986): instead of having the braid axis fixed and moving by isotopies the arcs of the link diagram for reaching  a closed  braid form, he considers the diagram fixed and instead he `threads' the braid axis (in the form of a simple closed curve projected on the  plane of the link diagram) through the arcs of the diagram, so that each subarc winds around the axis in the counterclockwise sense. In the same paper Morton uses a more technical version of his braiding algorithm for proving the Markov theorem.

A novel approach to the  Alexander theorem is due to Shuji Yamada  (\cite{Ya}, 1987) using the auxilliary concept of the `Seifert circles'. The  Seifert circles of a link diagram are created after smoothing all crossings according to the orientations of their arcs. Yamada introduces grouping operations on the system of Seifert circles which correspond to isotopies on the link diagram and which terminate with a braided form. 

\begin{wrapfigure}{R}{0.60\textwidth} 
\centering 
\includegraphics[width=0.57\textwidth]{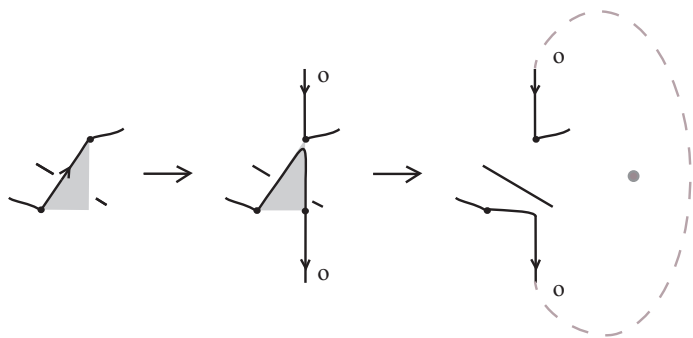} 
  \caption{An $L$-braiding move results in a pair of corresponding braid strands}
   \label{basicbraiding} 
\end{wrapfigure}

An additional value of Yamada's braiding algorithm is that it implies  equality between the Seifert number and the braid index of the link.  Pierre Vogel gave in 1990  a more sophisticated braiding algorithm based on the one by Yamada, where trees are used for measuring complexity  \cite{Vo}. Vogel's algorithm was then used by Pawel Traczyk for proving the Markov theorem (\cite{Tr}, 1992). A further approach to the Alexander theorem was made by the author in the 1990's \cite{La1,La2,LR1}. This algorithm results in open braids that one can `read' directly.  Here  we mark the local maxima and minima of the link diagram with respect to the height function and the opposite subarcs are now called {\it up-arcs}, so that we can have for each up-arc an orthogonal sliding triangle of the same type, `under' or `over', as the up-arc. 

\begin{wrapfigure}{L}{0.6\textwidth}  
\centering 
\includegraphics[width=0.57\textwidth]{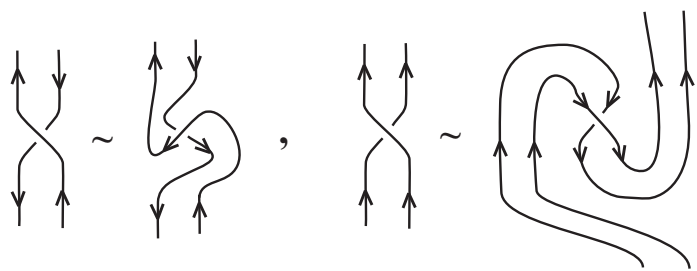} 
  \caption{A half twist and a full twist on crossings}
   \label{twistxings} 
\end{wrapfigure}

\begin{wrapfigure}{R}{0.7\textwidth} 
\centering 
\includegraphics[width=0.65\textwidth]{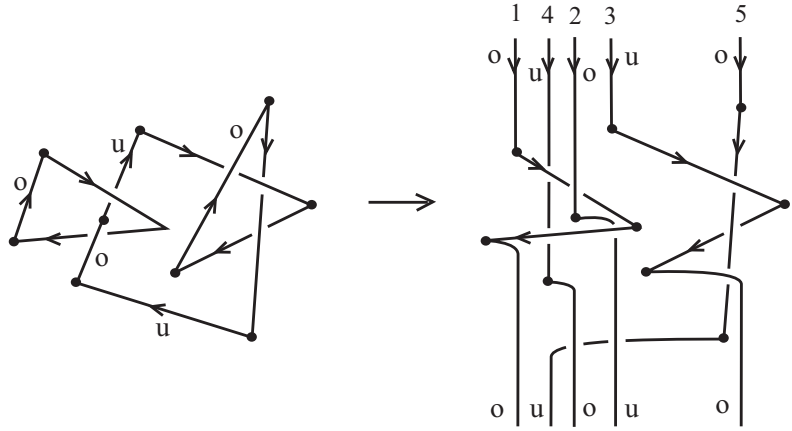} 
\caption{An example of applying the $L$-braiding algorithm}
\label{breg} 
\end{wrapfigure}

 The elimination of an up-arc is an {\it $L$-braiding move} and it consists in cutting the arc at some point, say the upmost, and then pulling the two ends, the upper upward and the lower downward, and keeping them aligned,  so as to obtain a pair of corresponding braid strands, both running entirely {\it over\/}
the rest of the diagram  or entirely {\it under\/} it, according to the type of the up-arc.   Fig.~\ref{basicbraiding} illustrates the abstraction of an `over'   $L$-braiding move and Fig.~\ref{breg} an example of applying the $L$-braiding algorithm. The closure of the resulting tangle is a link diagram, obviously isotopic to the original one, since from the up-arc we created a stretched loop  isotopic to the original up-arc. Any type of closure can apply here too.

 However, if we want to apply the braiding algorithm to the closure of a braid and obtain  the initial braid back, it is convenient to apply the vertical closure. 
 Finally, a really elegant version of the above braiding algorithm is given by the author and Louis H. Kauffman (\cite{KL1}, 2004): it uses the basic $L$-braiding move described above for up-arcs with no crossings, but it braids the up-arcs with crossings simply by rotating them on their projection plane, as illustrated in Fig.~\ref{twistxings}. 

Each algorithm for proving the  Alexander theorem has its own flavour and its own advantages and disadvantages but they  are all based on the same idea. Namely, to eliminate one by one the arcs of the diagram that have the wrong sense with respect to the chosen braid axis and to replace them by admissible ones.

\section{The braid group}

\begin{wrapfigure}{R}{0.65\textwidth} 
\centering 
\includegraphics[width=0.6\textwidth]{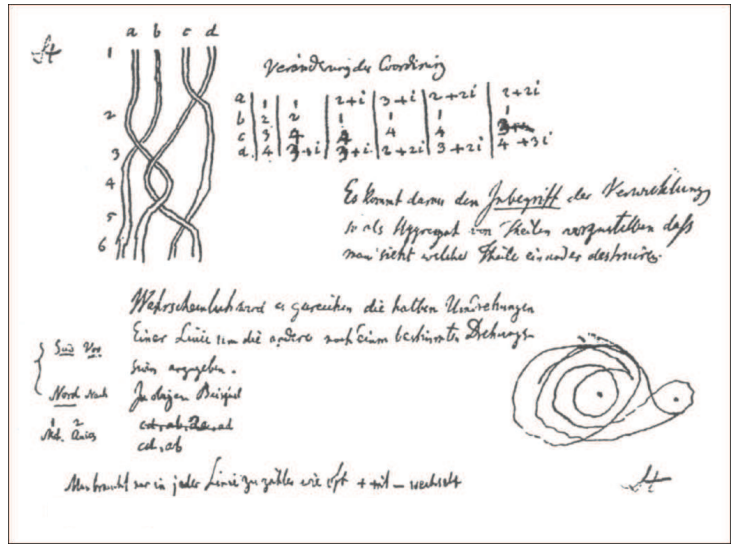} 
\caption{Gauss' handwritten study of braids} 
\label{gauss}
\end{wrapfigure}

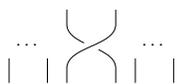
\begin{wrapfigure}{L}{0.20\textwidth} 
\scalebox{0.65}{
\begin{tikzpicture}[scale=1]
\braid[number of strands =1] 1  ;
\end{tikzpicture}
\begin{tikzpicture}[scale=1]
\node [label={[shift={(0.0, 0.3)}]$\cdots$}] {};
\end{tikzpicture}
\begin{tikzpicture}[scale=1]
\braid[number of strands =1] 1  ;
\end{tikzpicture}
\quad
\begin{tikzpicture}[scale=1]
\braid s_1^{-1}   ; 
\end{tikzpicture}
\quad
\begin{tikzpicture}[scale=1]
\braid[number of strands =1] 1  ;
\end{tikzpicture}
\begin{tikzpicture}[scale=1]
\node [label={[shift={(0.0, 0.3)}]$\cdots$}] {};
\end{tikzpicture}
\begin{tikzpicture}[scale=1]
\braid[number of strands =1] 1  ;
\end{tikzpicture}
}
\caption{The braid generator $\sigma_i$}
   \label{brgenrs} 
\end{wrapfigure}

The systematic theory of braids was introduced by E. Artin  in 1925 \cite{Ar1,Ar2}. He started with the geometric definition and then studied the group structure of the set of equivalence classes of braids on $n$ strands, after analyzing the braid isotopies with respect to the height function. These equivalence classes are also called {\it braids} and they form the {\it classical braid group $B_n$} with  operation the concatenation.  The group $B_n$ has a presentation of the  group $B_n$ with finitely many  generators $\sigma_i$, $i=1,\ldots,n-1$, and a finite set of defining relations [Artin, Chow]:
\[
\sigma_{i}\sigma_{i+1}\sigma_{i}= \sigma_{i+1}\sigma_{i}\sigma_{i+1}\qquad \text{and}\qquad
\sigma_i\sigma_j = \sigma_j\sigma_i \quad \text{for}\quad \vert i-j\vert >1.
\]
  The generators $\sigma_i$ resemble the elementary transpositions $s_i$ of the symmetric group $S_n$. They can be viewed as elementary braids of one crossing between the consecutive strands $i$ and $i+1$, also carrying   the topological information of which strand crosses over the other, see Fig.~\ref{brgenrs}. The  most important braid relation is the braided Reidemeister III move (all arrows down), while the Reidemeister II move is a direct consequence of the fact that the generators $\sigma_i$ are invertible. 

The group $B_n$ surjects on the symmetric group $S_n$ by corresponding the generator $\sigma_i$ to the elementary transposition $s_i$. The group $S_n$ has the extra relations $s_i^2 = 1$, $i=1,\ldots,n-1$, which are responsible for its finite order. The permutation induced by a braid tells us the number of link components of its closure by counting the disjoint cycles of the permutation. Braid groups play a central role in many areas of mathematics. A complete reference on braid groups and related topics can be found in the text by Kassel and Turaev \cite{KT}. 

It is worth noting that  C.F. Gauss was also thinking about the concept of a braid, probably in the frame of studying interwinding curves in space in the context of  his theory of electrodynamics.
 In Fig.~\ref{gauss} we see a handwritten note of Gauss, page~283 of his Handbuch~7, containing a sketch of a  braid that closes to a 3-component link, a coding table as well as a curve configuration winding around two points on its projection plane.

\section{Link isotopy and braid equivalence}\label{braidequiv}

We next consider equivalence relations in the set of all  braids  that correspond in the closures to link isotopy. This problem was first studied by A.A. Markov \cite{Ma}, after having available the Alexander theorem and  Artin's algebraic structure of the braid group. 

\begin{wrapfigure}{R}{0.65\textwidth} 
\centering 
\includegraphics[width=0.6\textwidth]{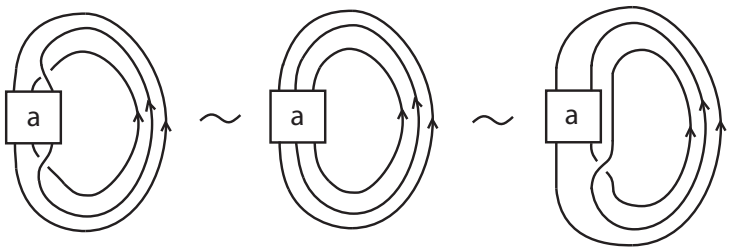} 
\caption{Closure of braid conjugation and braid stabilization induces isotopy} 
\label{conjugation}
\end{wrapfigure}

Closing a braid $a$ and a conjugate of $a$ by one braid generator, the two resulting links differ by Reidemeister II \& III moves that take place in the closure part of the conjugate of $a$. See Fig.~\ref{conjugation}. Similarly, adding a new strand at the end of the braid, which crosses once over or under the last strand corresponds in the closures to a kink, that is, to a Reidemeister I move. This  move is called `stabilization move'. See Fig.~\ref{conjugation}. Clearly, conjugations in all braid groups and stabilization moves are seemingly independent from each other and must figure in a braid equivalence that reflects link isotopy. The question is whether there are any other hidden moves for capturing the complexity of link isotopy. Yet, this is not the case, as we see from the following:

\begin{thm}[A.A. Markov, 1936] \label{markov} 
Two oriented links are isotopic if and only if any two corresponding braids differ by a finite sequence of braid relations and  the moves: 
\begin{center}
(i) \ Conjugation $\sigma_i^{-1} a \sigma_i \sim a$ and \ (ii) \ Stabilization $a \sigma_n^{\pm 1} \sim a$ for any $a \in B_n$ and for all $n\in {\mathbb N}$.
\end{center}
\end{thm}

The  statement of the theorem  included originally {\it three  moves}: the two local moves above and another more global one, the exchange move, that generalizes the Reidemeister II move. The sketch of proof by A.A. Markov \cite{Ma} used Alexander's braiding algorithm.
 Soon afterwards, N.~Weinberg reduced the exchange move to the {\it two moves} of the statement \cite{Wei}. The interest in the braid equivalence was rekindled by Joan Birman after following the talk of some unknown speaker. Birman produced a  rigorous proof of the theorem, filling in all details, by using a more technical version of Alexander's braiding algorithm \cite{Bi1}. A few years later Daniel Bennequin gave a different proof using 3-dimensional contact topology \cite{Ben}. In 1984 the Jones polynomial was discovered \cite{Jo1,Jo2}, a new powerful link invariant, whose construction used a representation of the braid group in the Temperley--Lieb algebra an the Alexander and Markov theorems. 
 This discovery led to new approaches to the Markov theorem. Hugh Morton gave a new proof using his threading algorithm \cite{Mo}, Pawel Traczyk proved the Markov theorem using Vogel's algorithm \cite{Tr}, and Joan Birman revisited the theorem with William Menasco  using Bennequin's ideas \cite{BM}.  Finally, in the 1990's, the author discovered a more geometric braid equivalence move, the $L$-move, and proved with Colin P. Rourke an {\it one-move} analogue of the Markov theorem, whose proof used the braiding moves described earlier \cite{La1, La2, LR1}:

\begin{thm} \label{L} 
Two oriented links are isotopic if and only if any two corresponding braids differ by a finite sequence of braid relations and the $L$-moves.
\end{thm}

\begin{wrapfigure}{L}{0.7\textwidth} 
\centering 
\includegraphics[width=0.62\textwidth]{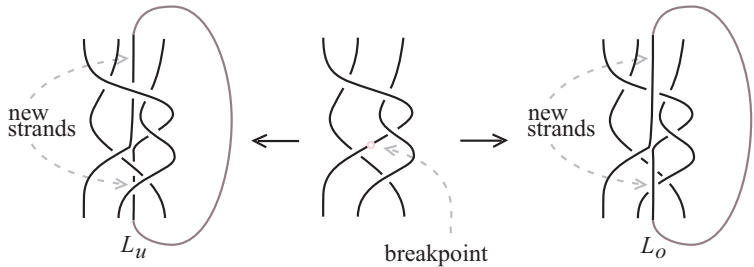} 
\caption{An $L_u$-move and an $L_o$-move at the same point} 
\label{Lmoves} 
\end{wrapfigure}

\begin{wrapfigure}{L}{0.65\textwidth} 
\centering
\includegraphics[width=0.6\textwidth]{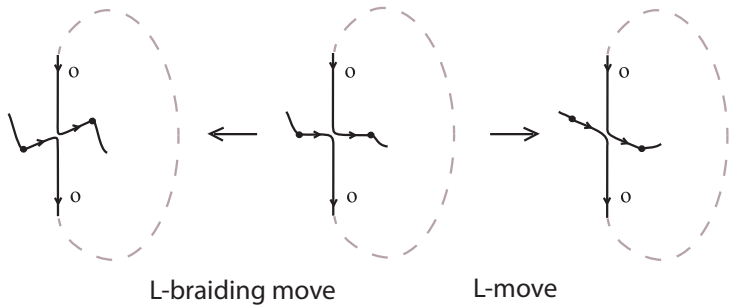}
\caption{The symmetry of the braiding and the $L$-move}
\label{symmetry}
\end{wrapfigure}

An {\it $L$-move} resembles an $L$-braiding move: it consists in cutting a braid arc at some point and then pulling the two ends, the upper downward  and the lower upward, keeping them aligned, and so as to obtain a new pair of corresponding braid strands, both running entirely {\it over\/}
the rest of the braid  or entirely {\it under\/} it, according to the type of the move, denoted $L_o$ or $L_u$ respectively.  Fig.~\ref{Lmoves} illustrates an example with both types of $L$-moves taking place at the same point of a braid. The closure of the resulting braid differs from the closure of the initial one by a stretched loop. View also Fig.~\ref{symmetry} for an abstact illustration of the similarity of the $L$-braiding move and the $L$-move. 

\begin{wrapfigure}{R}{0.65\textwidth} 
\centering 
\includegraphics[width=0.6\textwidth]{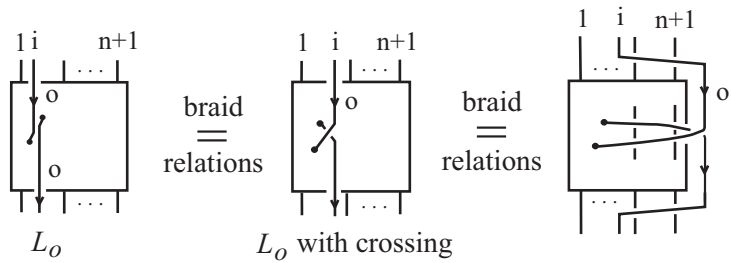} 
\caption{The $L$-moves have algebraic expressions} 
\label{algL}
\end{wrapfigure}

 The $L$-moves are geometric.  However, as we see in the middle illustration of  Fig.~\ref{algL}, using braid isotopy an $L$-move can be also viewed as introducing  a crossing inside the braid `box', so  in the closure it creates a stretched Reidemeister I kink. This way of viewing the $L$-moves  shows that they generalize the stabilization moves. It also renders them locality and it leads to the observation that they have algebraic expressions, as it is clear from Fig.~\ref{algL}.  Furthermore, it follows from Theorem~\ref{L}  that conjugation can be achieved by braid relations and  $L$-moves. 

\smallbreak

The Markov theorem  or Theorem~\ref{L}  are not easy to prove. For proving the `only if' part one needs first to take two diagrams of the two isotopic links and produce  corresponding  braids, using some braiding algorithm, and then show that the two braids are Markov equivalent (resp. $L$-move equivalent). In practice this means that any choices made on a given link diagram when applying the braiding algorithm correspond to Markov equivalent  (resp. $L$-move equivalent) braids and that if two link diagrams differ by an isotopy move the corresponding braids are also Markov equivalent (resp. $L$-move equivalent).  For this analysis it is crucial to have the isotopy moves local and the braiding moves independent from each other. In this way one can always assume to have done almost all braiding in the otherwise identical diagrams in question and to be left only with the isotopy move or algorithmic choice by which they differ. Then the two braid diagrams are directly comparable. 

The braiding algorithm of \cite{La1,La2,LR1} is particularly appropriate for proving the Markov theorem or its equivalent Theorem~\ref{L}, after enchancing it with some extra technicalities \cite{La1,La2,LR1} for ensuring independence of the sequence of the braiding moves. 
 This is because the $L$-braiding moves and the $L$-moves are simple and have a basic symmetric  interconnection, as illustrated in Fig.~\ref{symmetry}, so they comprise, in fact, one very  fundamental uniform move, enabling one to trace easily how the algorithmic choices and the isotopy moves on diagrams affect the final braids. 

% We shall be using the term `$L$-move' for both, the $L$-braiding move and the $L$-move.    

\section{Extensions to other diagrammatic settings}

Given a diagrammatic knot theory there are deep interrelations between the diagrammatic  isotopy in this theory, the braid structures and the corresponding braid equivalences.

More precisely, the isotopy moves that are allowed, but also moves that are forbidden in the setting, determine the corresponding braid isotopy, the way the closure of a braid is realized, and also the corresponding braid equivalence moves. 
 On the other hand, braid equivalence theorems are important for understanding the structure and classification of knots and links in various settings and for constructing invariants of knots and links using algebraic means, for example via Markov traces on quotient algebras of the corresponding  braid groups.

	The $L$-braiding  and the $L$-moves provide a uniform and flexible ground for formulating and proving braiding and braid equivalence theorems for any diagrammatic setting. Indeed, their simple and fundamental  nature together with the fact that the $L$-moves are geometric and can be localized  are the reasons that they can adapt to all diagrammatic categories where the notions of braid and diagrammatic isotopy are defined. This is particularly useful in settings where algebraic braid structures are not immediately apparent. 	Indeed, the statements are first geometric and then they gradually turn algebraic, if algebraic braid structures are available. 
 																																												\smallbreak
				
	The $L$-move techniques were first employed for proving braiding and braid equivalence theorems for classical knots in $3$-manifolds with or without boundary; namely in knot complements, in closed, connected, oriented (c.c.o.)  $3$-manifolds, which are obtained from the $3$-sphere, $S^3$, via the `surgery technique', as well as in handlebodies.

\begin{wrapfigure}{R}{0.72\textwidth} 
\centering
\includegraphics[width=.65\textwidth]{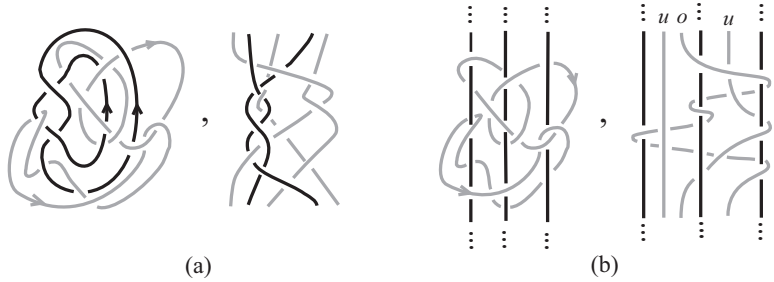}
\caption{(a) a mixed link and a geometric mixed braid for link complements and c.c.o. 3-manifolds; (b)  a mixed link and a geometric mixed braid for handlebodies}
\label{mfds}
\end{wrapfigure}

\begin{wrapfigure}{L}{0.65\textwidth}                                               
\centering
\includegraphics[width=0.6\textwidth]{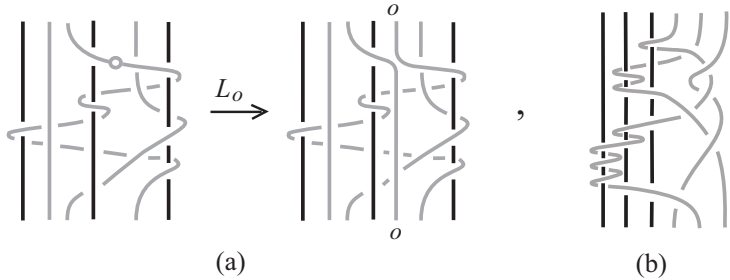}
\caption{An $L$-move in a mixed braid; (b) an algebraic mixed braid}
\label{Lmfds}
\end{wrapfigure}

  The idea here is to fix in $S^3$ a closed braid representation of the $3$-manifold  and then represent knots and braids in the $3$-manifold as {\it mixed links} and {\it mixed braids} in $S^3$, which contain the {\it fixed part}, representing the $3$-manifold, and the {\it moving part}, representing the knot/braid in the $3$-manifold. View Fig.~\ref{mfds} for concrete examples. Then, knot isotopy in the $3$-manifold is translated into mixed link isotopy in $S^3$, which applies only on the moving part.  In the case of c.c.o. $3$-manifolds we have the isotopy moves for the knot complements, as well as extra isotopy moves, the {\it band moves}, related to the surgery description of the manifold. 
			The mixed braid equivalence for knot complements comprises $L$-moves which take place only on the moving parts of mixed braids. See Fig.~\ref{Lmfds}(a), while in the case of c.c.o. $3$-manifolds we have the extra {\it braid band moves}.	Then the $L$-move braid equivalences turn into algebraic statements with the use of the  algebraic mixed braids \cite{La3}, see Fig.~\ref{Lmfds}(b) for an example. 																																											

 For the case of a handlebody we have the same setting as for a knot complement. The difference here is  that a knot may not pass beyond its boundary from either end, and this is reflected both in the definition of the closure of a mixed braid as well as in the corresponding braid equivalence. Namely, the closure is realized by simple closing arcs, slightly tilted in the extremes, which run `over' or `under' the rest of the diagram, and different choices may lead to non-isotopic closures. Furthermore, in the mixed braid equivalence some `loop' conjugations are not allowed.  Details on the above can be found in \cite{LR1,La3,LR2,DL,HL,La4,La5}. See also \cite{Su} and \cite{Sk}.

	\smallbreak
	
\begin{wrapfigure}{R}{0.55\textwidth} 
\centering
\includegraphics[width=0.5\textwidth]{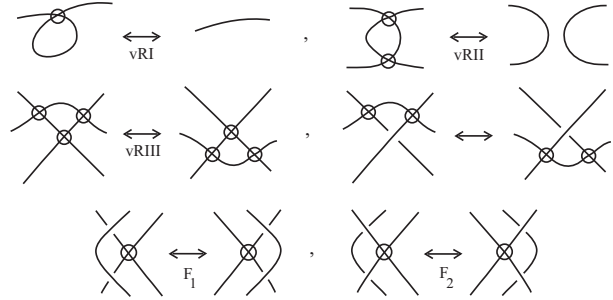}
\caption{Virtual moves: allowed and forbidden}
\label{vmoves}
\end{wrapfigure}

The next application of the $L$-move methods was in virtual knot theory.  Virtual knot theory was introduced by Louis~H. Kauffman \cite{Kau2} and it is an extension of classical  knot theory. In this extension one  adds  a `virtual' crossing that is neither an over-crossing nor an under-crossing. Fig.~\ref{vmoves} illustrates the diagrammatic moves that contain virtual crossings. In this theory we have the {\it virtual forbidden moves}, F1 and F2, with two real crossings and one virtual. We have here the virtual braid group \cite{Kau2,KL1,Bar}, which extends the classical braid group.	The forbidden moves make it harder to braid a virtual knot diagram, so the idea of rotating crossings that contain up-arcs before braiding a diagram (recall Figure~\ref{twistxings}) comes in handy in this setting. The interpretation of an $L$-move as introducing an in-box crossing proved very crucial  in the search of the types of $L$-moves needed in the setting, as they are related to the types of kinks allowed in the given isotopy. 

 So, we have $L$-moves introducing a real or a virtual crossing facing to the right or to the left of the braid. Moreover, the presence of the forbidden moves in the theory leads to the requirement that the strands of an $L$-move cross the other strands of the virtual braid only virtually, and also to a type of virtual $L$-move coming from a `trapped' virtual kink.  The above led in \cite{KL2}  to formulations of virtual braid equivalence theorems for the virtual braid group, for the welded braid group \cite{FRR} and for some analogues of these structures, complementing the prior results of Seiichi Kamada in \cite{Ka}, where the more global exchange moves are used.

\begin{wrapfigure}{L}{0.77\textwidth} 
\centering
\includegraphics[width=0.73\textwidth]{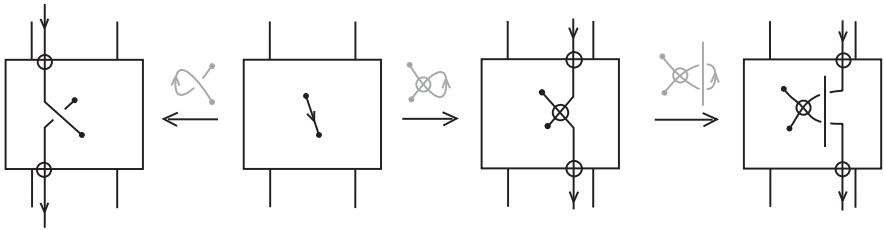}
\caption{Types of virtual $L$-moves}
\label{virtualL}
\end{wrapfigure}

	Furthermore, the $L$-move techniques have been used  by Vassily Manturov and Hang Wang for formulating a Markov-type theorem for free links \cite{mawa}. Also, by Carmen Caprau and co-authors for obtaining a braid equivalence for virtual singular braids
 \cite{capega} as well as for virtual trivalent braids \cite{cadiposa,cacoda}. 

		\smallbreak
		
\begin{wrapfigure}{R}{0.5\textwidth} 
\centering
\includegraphics[width=0.45\textwidth]{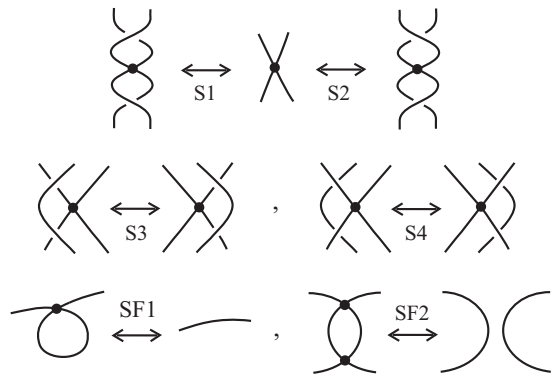}
\caption{Singular moves: allowed and forbidden}
\label{smoves}
\end{wrapfigure}

Singular knot theory is related to Vassiliev's theory of knot invariants. Fig.~\ref{smoves} illustrates the diagrammatic moves in the theory as well as the {\it singular forbidden moves}, SF1 and SF2. The singular crossings together with the real crossings and their inverses generate the `singular braid monoid' introduced in different contexts by Baez\cite{Ba}, Birman\cite{Bi2} and Smolin\cite{Sm}. Braiding a singular knot diagram becomes particularly simple by using the idea of rotating singular crossings that contain up-arcs before braiding (Figure~\ref{twistxings}) and the $L$-braiding moves. 
  An algebraic singular braid equivalence is proved by Bernd Gemein in \cite{Ge} and, assuming this result, in \cite{La4} the $L$-move analogue is formulated. Clearly, there is no $L$-move introducing a singular crossing, as the closure of such a move would contract to a kink with a singular crossing, and this is not an isotopy move in the theory. Also, there is no conjugation by a singular crossing, since this is not an invertible element in the monoid; yet, we can talk about `commuting' in the singular braid monoid: $ab \sim ba$. 
	
	\smallbreak

\begin{wrapfigure}{L}{0.6\textwidth} 
\centering
\includegraphics[width=0.55\textwidth]{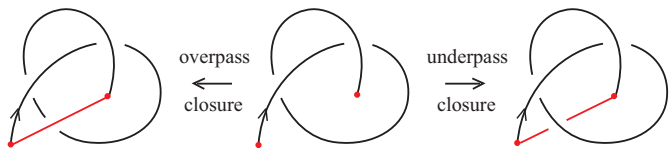}
\caption{A knotoid and its two closures to knots}
\label{knotoid}
\end{wrapfigure}
	
\begin{wrapfigure}{R}{0.45\textwidth} 
\centering 
\includegraphics[width=0.4\textwidth]{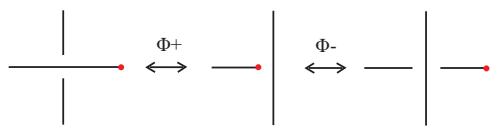} 
\caption{The forbidden moves for knotoids} 
\label{forbidoid} 
\end{wrapfigure}

Another very interesting diagrammatic category is the theory of knotoids and braidoids.
	The theory of knotoids was introduced by Vladimir Turaev 
 in 2012 \cite{Tu}. A knotoid diagram
 is an open curve in an oriented surface, having finitely many self-intersections that are endowed with under/over data and with its two endpoints possibly lying in different regions of the diagram. For an example see middle illustration of Fig.~\ref{knotoid}. The theory of knotoids is a complex diagrammatic theory, and its complexity  lies in the {\it knotoid forbidden moves}, $\Phi_+$ and $\Phi_-$, that prevent the endpoints from slipping under or  over other arcs of the diagram. See Fig.~\ref{forbidoid}. 

The theory of  spherical knotoids (i.e., knotoids in the two-sphere) extends the theory of  classical knots and also proposes a new diagrammatic approach to classical knots, which arise via the `overpass' or `underpass' closures, see Fig.~\ref{knotoid}. This approach promises reducing of the computational complexity of knot invariants  \cite{Tu}. On the other hand, planar knotoids surject to spherical knotoids, but do not inject. This means that planar knotoids provide a much richer combinatorial structure than the spherical ones. This fact has interesting implications in the study of proteins \cite{GDBS,GGLDSK}.

\begin{wrapfigure}{L}{0.4\textwidth} 
\centering 
\includegraphics[width=0.35\textwidth]{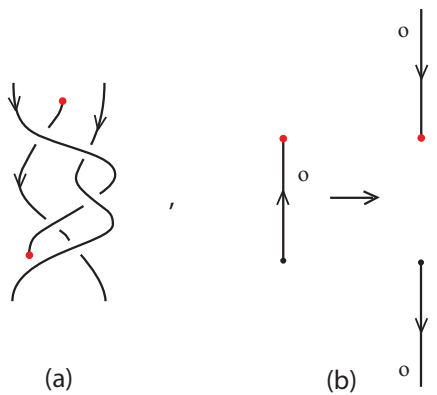} 
\caption{(a) A braidoid; (b) $L$-braidoiding at an endpoint} 
\label{braidoids} 
\end{wrapfigure}

Recently, the theory of braidoids has been introduced and developed \cite{GL1,GL2}, which  extends the classical braid theory.  A `braidoid' is like a classical braid, but two of its strands terminate at the endpoints. For an example see Fig.~\ref{braidoids}(a). The forbidden moves play a role in the algorithm for turning planar knotoids to braidoids and they affect the definition of the closure operation on braidoids, in which the endpoints do not participate. Namely, we close corresponding braidoid ends using vertical arcs with slightly tilted extremes, running over or under the rest of the diagram, and this needs to be specified in advance, as different choices may result in non-isotopic knotoids (due to the forbidded moves). For turning a planar knotoid into a braidoid, we use the $L$-braiding moves for up-arcs not containing an endpoint and the analogous moves illustrated in Fig.~\ref{braidoids}(b) (with choice `o' in the figure) for the ones that contain an endpoint. 
For  a braidoid equivalence we use the $L$-moves that can take place at any point of a braidoid except for the endpoints. We note that for the braidoids we do not have yet an appropriate algebraic structure, so the $L$-equivalence is as good as we can have so far.

It is worth adding at this point that, in \cite{GK1} Neslihan G\"ug\"umc\"u and Louis Kauffman give a faithful lifting of a planar knotoid to a space curve, such that the endpoints remain attached to two paralell lines. In \cite{KoLa1} Dimitrios Kodokostas and the author make the observation that 
 this interpretation of planar knotoids is related to the knot theory of the handlebody of genus two and these ideas are further explored in \cite{KoLa2}, where the notion of `rail knotoid' is introduced.

Further, in \cite{La7} the $L$-move techniques are applied to long knots. 

\smallbreak

Finally, in \cite{Sch} Nancy Scherich provides a computer
implemented, grid diagrammatic proof of the Alexander theorem, based on the $L$-braiding moves. Another result of analogous flavour is a Markov-type theorem for ribbon torus-links in ${\mathbb R}^4$ by Celeste Damiani \cite{Da}.

	\smallbreak
	
Surveys on many of the above results are included in \cite{La4,La5,GKL}, while a more complete presentation is to appear \cite{La6}.

%%%%%%%%%%%%%%%%%%%%%%%%%%%%%%%%%%%%%%%%%%%%%%%%%%%%%%%%%%%%%%%%%%%%%%%%%%%%%%%%%%%%%

\end{document}